\documentclass[reqno]{amsart}
\usepackage{amsmath}
\usepackage{xspace}
\usepackage{alltt}

\theoremstyle{definition}
\newtheorem{definition}{Definition}
\theoremstyle{plain}
\newtheorem{thm}{Theorem}

\numberwithin{equation}{section}
\allowdisplaybreaks[1]

\newcommand{\qBin}[3]{\genfrac{[}{]}{0pt}{0}{#1}{#2}_{#3}}
\newcommand{\qBinsm}[3]{\genfrac{[}{]}{0pt}{1}{#1}{#2}_{#3}}

\newcommand{\myvec}[1]{\mathbf{#1}}
\newcommand{\ivec}{\myvec{i}}
\newcommand{\jvec}{\myvec{j}}
\newcommand{\Ivec}{\myvec{I}}
\newcommand{\Jvec}{\myvec{J}}
\newcommand{\kvec}{\myvec{k}}
\newcommand{\nvec}{\myvec{n}}

\newcommand{\kfree}{$\kvec$\nobreakdash-free\xspace}
\newcommand{\qhg}{$q$\nobreakdash-hyper\-geometric\xspace}

\newcommand{\lbr}{\symbol{123}}
\newcommand{\rbr}{\symbol{125}}
\newcommand{\mysub}[1]{\(\sb{\mathtt{#1}}\)}
\newlength{\axellength}
\newcommand{\In}[1]{\settowidth{\axellength}{\texttt{12345#1}}%
\makebox[\axellength][r]{\textsf{\scriptsize In[#1]:=}}}
\newcommand{\Out}[1]{\settowidth{\axellength}{\texttt{12345#1}}%
\makebox[\axellength][r]{\textsf{\scriptsize Out[#1]=}}}

\DeclareMathOperator{\SUM}{SUM}

\def\eqn#1{\eqref{eq:#1}}

\newenvironment{myhack}%
  {\renewcommand{\normalfont}%
    {\usefont\encodingdefault\familydefault\seriesdefault\shapedefault\normalsize}}%
  {\renewcommand{\normalfont}%
    {\usefont\encodingdefault\familydefault\seriesdefault\shapedefault}}

\begin{document}

\title[A Polynomial Identity Implying a Partition Theorem of G\"ollnitz]%
      {A Computer Proof of a Polynomial Identity\\
       Implying a Partition Theorem of G\"ollnitz}
\author[A.~Berkovich]{Alexander Berkovich}
\address{Department of Mathematics, The Pennsylvania State University,
         University Park, PA~16802, USA}
\email{alexb@math.psu.edu}
\thanks{The first author was partially supported by SFB-grant
        F1305 of the Austrian FWF}
\author[A.~Riese]{Axel Riese}
\address{Research Institute for Symbolic Computation,
         Johannes Kepler University,
	 A--4040 Linz, Austria}
\email{Axel.Riese@risc.uni-linz.ac.at}
\thanks{The second author was supported by SFB-grant
        F1305 of the Austrian FWF}

\begin{abstract}
In this paper we give a computer proof of a new polynomial identity, which extends a
recent result of Alladi and the first author. In addition, we provide computer
proofs for new finite analogs of Jacobi and Euler formulas. All computer proofs
are done with the aid of the new computer algebra package \texttt{qMultiSum} developed by
the second author. \texttt{qMultiSum} implements an algorithmic refinement of Wilf
and Zeilberger's multi-$q$-extension of Sister Celine's technique utilizing
additional ideas of Verbaeten and Wegschaider.
\end{abstract}

\maketitle

\section{G\"ollnitz's Partition Theorem and Related\\ $q$-Hypergeometric Identities}

In 1967, G\"ollnitz~\cite{G} proved the following deep partition theorem:

\begin{thm} \label{thm:1}
Let $A(N)$ denote the number of partitions of $N$ in the form $N=n_1+n_2+n_3+\dotsb$,
such that no part is equal to $1$ or $3$, and $n_i-n_{i+1}\ge 6$ with strict inequality if
$n_i\equiv 0,1,3 \pmod 6$.\\
Let $B(N)$ denote the number of partitions of $N$ into distinct parts $\equiv 2,4,5 \pmod 6$.
Then
\[
   A(N)=B(N).
\]
\end{thm}

In \cite{AAG}, Alladi, Andrews, and Gordon reformulated and refined Theorem~\ref{thm:1} using
the language of colored partitions.
To state their theorem, we will need a few definitions.

Let $P(N,i,j,k)$ denote the number of partitions of $N$ into parts occurring in three
(primary) colors ordered as
\begin{equation} \label{eq:1.1}
   \mathbf{A}<\mathbf{B}<\mathbf{C},
\end{equation}
such that parts in the same color are distinct and the number of parts in colors
$\mathbf{A},\mathbf{B},\mathbf{C}$ is equal to $i,j,k$, respectively.

Next, consider partitions $\pi$, such that parts equal to $1$ may occur in three primary
colors~\eqn{1.1}, but parts $\ge2$ may occur in six colors ordered as
\begin{equation} \label{eq:1.2}
   \mathbf{AB}<\mathbf{AC}<\mathbf{A}<\mathbf{BC}<\mathbf{B}<\mathbf{C}.
\end{equation}
In addition, the gap between adjacent parts is required to be $\ge1$, where the gap
may equal $1$ only if both parts are either of the
same primary color or the larger part is in a color of higher order
according to~\eqn{1.2}.
These partitions $\pi$ were called Type--$1$ partitions in~\cite{AAG}.

We can now state the following result~\cite{AAG}.

\begin{thm} \label{thm:2}
Let $G(N,a,b,c,ab,ac,bc)$ denote the number of Type--$1$ partitions
of $N$ with $a$ parts in color $\mathbf{A}$, $\ldots$, $bc$ parts in color
$\mathbf{BC}$.
\textup{(}Note that $bc$ is not $b\cdot c$!\textup{)}
Then
\[
\sum_{\substack{i,j,k \\ \text{constraints}}} G(N,a,b,c,ab,ac,bc)=P(N,i,j,k),
\]
where the sum on the left is over the variables $a,b,c,ab,ac,bc$ subject to the
$i,j,k$-constraints,
which here \textup{(}and throughout\textup{)} are
\[
   \begin{cases} i=a+ab+ac,\\ j=b+ab+bc,\\ k=c+ac+bc.\end{cases}
\]
\end{thm}

To see the connection between Theorem~\ref{thm:1} and Theorem~\ref{thm:2},
we denote part $n$ in color $\mathbf{A}$
as $\mathbf{A}_n$, $\ldots$, part $n$ in color $\mathbf{BC}$ as $\mathbf{BC}_n$.
Next, we replace the colored integers $\mathbf{A}_n,\ldots,\mathbf{BC}_n$
by regular integers according to the following rules:
\begin{equation} \label{eq:1.3}
   \left\{ \hspace*{-\arraycolsep} \begin{array}{rl} 
   \mathbf{A}_n \rightarrow 6n-4,\\
   \mathbf{B}_n \rightarrow 6n-2,\\
   \mathbf{C}_n \rightarrow 6n-1,\\[1ex]
   \mathbf{AB}_n \rightarrow 6n-6, & n>1,\\
   \mathbf{AC}_n \rightarrow 6n-5, & n>1,\\
   \mathbf{BC}_n \rightarrow 6n-3, & n>1.\end{array}\right.
\end{equation}

Note that this replacement converts the ordering
\[ 
  \mathbf{AB}_n<\mathbf{AC}_n<\mathbf{A}_n<\mathbf{BC}_n<\mathbf{B}_n<\mathbf{C}_n,
  \quad n>1
\]
into the conventional ordering
\[ 
  6n-6<6n-5<6n-4<6n-3<6n-2<6n-1, \quad n>1.
\]
Recalling that part $1$ can occur only as $\mathbf{A}_1,\mathbf{B}_1,\mathbf{C}_1$, we infer
that no conventional part equals $1$ or $3$.
Also, one can easily check that under the transformations~\eqn{1.3} the color-gap
conditions become identical with the gap conditions in Theorem~\ref{thm:1}.
So, summing  over $i,j,k$ one immediately obtains the G\"ollnitz partition theorem.

\medskip
To prove Theorem~\ref{thm:2}, the authors of~\cite{AAG} stated it in an analytic form,
which they called \textit{Key Identity}, as follows:
\begin{equation} \label{eq:1.7}
   \sum_{\substack{i,j,k\\\text{constraints}}}
   \frac{q^{T_t+T_{ab}+T_{ac}+T_{bc-1}}
   (1-q^a+q^{a+bc})} {(q)_a(q)_b(q)_c(q)_{ab}(q)_{ac}(q)_{bc}}=
   \frac{q^{T_i+T_j+T_k}}{(q)_i(q)_j(q)_k},
\end{equation}
where
\begin{gather*} 
   t=a+b+c+ab+ac+bc,\\
   T_m=\frac{m(m+1)}{2},
\end{gather*}
and for $n\in\mathbb{Z}$, the $q$-shifted factorial of $a$ is defined as
\[ 
   (a)_n = (a;q)_n= \begin{cases}
   \prod^{n-1}_{j=0}(1-aq^j), & \text{if } n>0,\\
   1, & \text{if } n=0,\\
   \prod^{-n}_{j=1} {(1-aq^{-j})}^{-1}, & \text{if } n<0.\end{cases}
\]

Several proofs of the Key Identity~\eqn{1.7} have appeared in the literature
\cite{AA, AAG, AB1, R}.
In~\cite{AA}, Alladi and Andrews provided a straightforward
$q$-hypergeometric proof of~\eqn{1.7}.
Their proof made essential use of Jackson's
$q$-analog of Dougall's summation formula. (See, for instance, Gasper
and Rahman~\cite[(II.21), p.~238]{GR}.)
In \cite{R}, Riese used his computer algebra
package \texttt{qMultiSum} to find a very simple recursive proof of~\eqn{1.7}.
In~\cite{AB1}, Alladi and Berkovich proposed and proved a double bounded polynomial
generalization of~\eqn{1.7}:

\begin{thm} \label{thm:3}
If $i,j,k,L,M$ are integers, then
\begin{align}\label{eq:1.11}
   &\sum_{\substack{i,j,k\\ \text{constraints}}}
   q^{T_t+T_{ab}+T_{ac}+T_{bc-1}} \notag\\
   & \quad \times \Bigg\{ q^{bc} \qBin{L-t+a}{a}{q} \qBin{L-t+b}{b}{q}
   \qBin{M-t+c}{c}{q} \qBin{L-t}{ab}{q} \qBin{M-t}{ac}{q}
   \qBin{M-t}{bc}{q} \notag\\
   & \qquad \quad + \qBin{L-t+a-1}{a-1}{q} \qBin{L-t+b}{b}{q}
   \qBin{M-t+c}{c}{q} \qBin{L-t}{ab}{q} \qBin{M-t}{ac}{q}
   \qBin{M-t}{bc-1}{q} \Bigg\} \notag\\
   & = \sum_{s\ge0} q^{s(M+2)-T_s+T_{i-s}+T_{j-s}+T_{k-s}}
   \qBin{L-s}{s,i-s,j-s}{q} \qBin{M-i-j}{k-s}{q}.
\end{align}
\end{thm}

In~\eqn{1.11} and in what follows, the $q$-binomial and $q$-multinomial coefficients are
defined by
\[ 
  \qBin{n+m}{n}{q} = \begin{cases}
  \frac{(q^{m+1})_n}{(q)_n}, & \text{if } n\ge0,\\
   0, & \text{otherwise},\end{cases}
\]
and
\[ 
  \qBin{L}{a_1,a_2,\dots,a_n}{q} =
  \qBin{L}{a_1}{q} \qBin{L-a_1}{a_2}{q} \cdots \qBin{L-a_1-a_2-\dots-a_{n-1}}{a_n}{q}.
\]

The proof of Theorem~\ref{thm:3} given by Alladi and Berkovich~\cite{AB1} used recurrences
along with Jackson's $q$-Dougall formula.
It is easy to check that in the limit $L,M\to\infty$,
\eqn{1.11} reduces to~\eqn{1.7}.
Furthermore, if $L=M$, then the sum on the right
in~\eqn{1.11} can be evaluated with the aid of the $q$-Pfaff-Saalsch\"utz formula.
(See, for instance, Gasper and Rahman~\cite[(II.12), p.~237]{GR}.)
As a result, \eqn{1.11} becomes
\begin{align} \label{eq:1.14}
   &\sum_{\substack{i,j,k\\ \text{constraints}}}
   q^{T_t+T_{ab}+T_{ac}+T_{bc-1}} \notag\\
   & \quad \times \Bigg\{ q^{bc} \qBin{L-t+a}{a}{q} \qBin{L-t+b}{b}{q}
   \qBin{L-t+c}{c}{q} \qBin{L-t}{ab}{q} \qBin{L-t}{ac}{q}
   \qBin{L-t}{bc}{q} \notag\\
   & \qquad \quad + \qBin{L-t+a-1}{a-1}{q} \qBin{L-t+b}{b}{q}
   \qBin{L-t+c}{c}{q} \qBin{L-t}{ab}{q} \qBin{L-t}{ac}{q}
   \qBin{L-t}{bc-1}{q} \Bigg\} \notag\\
   & = q^{T_i+T_j+T_k} \qBin{L-k}{i}{q} \qBin{L-i}{j}{q} \qBin{L-j}{k}{q}.
\end{align}

A partition theoretical interpretation of~\eqn{1.14}, given in~\cite{AB1},
can be stated as follows:
\begin{thm} \label{thm:4}
Let $G_L(N,a,b,c,ab,ac,bc)$ denote $G(N,a,b,c,ab,ac,bc)$
with the additional constraint that no part exceeds $\mathbf{C}_L$.
Let $P_L(N,i,j,k)$ denote $P(N,i,j,k)$ with the additional constraints
\[
  \lambda(\mathbf{A})\le\mathbf{A}_{L-k},\quad
  \lambda(\mathbf{B})\le\mathbf{B}_{L-i},\quad
  \lambda(\mathbf{C})\le\mathbf{C}_{L-j},
\]
where $\lambda(\mathbf{A})$ is the largest part in color $\mathbf{A}$, and
$\lambda(\mathbf{B}),~\lambda(\mathbf{C})$ have the analogous interpretation.
Then, for $L\ge \max(i+j,j+k,k+i)$,
\[
   \sum_{\substack{i,j,k \\\text{constraints}}} G_L(N,a,b,c,ab,ac,bc)=P_L(N,i,j,k).
\]
\end{thm}

It was pointed out in~\cite{AB1} that \eqn{1.14} and Theorem~\ref{thm:4} can be employed to
derive new finite versions of many classical $q$-hypergeometric identities,
including those of Gauss, Jacobi and Lebesgue.

Since Theorem~\ref{thm:2} deals with partitions into parts occurring in three
primary colors, one may
suspect that a polynomial analog of~\eqn{1.7} should depend on three finitization
parameters $L_1,L_2,M$ and that \eqn{1.11} is just the special case $L_1=L_2=L$ of
this more general result.
And, indeed, further investigations led us to the following triple bounded
polynomial generalization of~\eqn{1.7}:

\begin{thm} \label{thm:5}
Let
\begin{align} \label{eq:1.16}
   &g_{i,j,k}(L_1,L_2,M) :=  \sum_{\substack{i,j,k\\\text{constraints}}}
   q^{T_t+T_{ab}+T_{ac}+T_{bc-1}}\notag\\
   & \times \Bigg\{ q^{bc} \qBin{L_1-t+a}{a}{q} \qBin{L_2-t+b}{b}{q}
   \qBin{L_2-t}{ab}{q} \qBin{M-t+c}{c}{q} \qBin{M-t}{ac}{q}
   \qBin{M-t}{bc}{q} \notag\\
   & \qquad + \qBin{L_1-t+a-1}{a-1}{q} \qBin{L_2-t+b}{b}{q}
   \qBin{L_2-t}{ab}{q} \qBin{M-t+c}{c}{q} \qBin{M-t}{ac}{q}
   \qBin{M-t}{bc-1}{q} \Bigg\}
\end{align}
and
\begin{align} \label{eq:1.17}
   & p_{i,j,k}(L_1,L_2,M) :=\sum_{s\ge0}q^{s(M+2)-T_s+T_{i-s}+T_{j-s}+T_{k-s}} \notag\\
   & \quad \times \qBin{L_1-s}{i-s}{q} \qBin{L_2-i}{j-s}{q}
   \qBin{L_2-i-j+s}{s}{q} \qBin{M-i-j}{k-s}{q}.
\end{align}
Then
\begin{equation} \label{eq:1.15}
g_{i,j,k}(L_1,L_2,M)=p_{i,j,k}(L_1,L_2,M).
\end{equation}
\end{thm}

We wish to comment that if any one of the parameters $i,j,k$ is set to zero,
then identity~\eqn{1.15} reduces to the double bounded key identity for
Schur's partition theorem~\cite{AB2}.

\medskip
While it is straightforward to extend the analysis of~\cite{AB1} to prove~\eqn{1.15},
our goal here is different.
We would like to use \eqn{1.15} as the testing ground for the package
\texttt{qMultiSum}, which has been recently developed by Riese~\cite{R}.
This package is described in some detail in Section~\ref{sec:qMultiSum}.
In Section~\ref{sec:triple}, we will use \texttt{qMultiSum} to obtain
a nested recursive proof of Theorem~\ref{thm:5}.
In Section~\ref{sec:Euler}, we will give a computer proof
of the new finite version of Jacobi's formula~\cite{AB1} and then, propose and prove
(automatically) a new finite analog of Euler's formula.
In Section~\ref{sec:Conclusion} we will conclude with a brief
discussion of the ``human insight", which went
into the computer proof of~\eqn{1.15}, and with a short description of
the prospects for future work.

\section{The Package \texttt{\upshape qMultiSum}} \label{sec:qMultiSum}

The object of this section is to give a short account on the
\textsf{Mathematica} package \texttt{qMultiSum}\footnote{available at
\textit{http://www.risc.uni-linz.ac.at/research/combinat/risc/software/qMultiSum}}
which computes recurrences for \qhg multi-sums.
The package has been written by the second author and will be described
only briefly here.
For more details the reader is referred to the forthcoming
article~\cite{R}.

The implementation is based on the method of \kfree recurrences, also known
as Sister Celine's technique (developed by Wilf and Zeilberger~\cite{WilfZeilberger:AlgoProof}).
For reasons of efficiency we also incorporated ideas  from
Wegschaider's~\cite{Wegschaider:Dipl} package
\texttt{MultiSum}\footnote{available at
\textit{http://www.risc.uni-linz.ac.at/research/combinat/risc/software/MultiSum}}
for ordinary hypergeometric summation.

\medskip
Let $\nvec = (n_1, \dots, n_s)$ and $\kvec = (k_1, \dots, k_r)$ be vectors
of variables ranging over the integers.
The central concept of (the $q$-version of) Sister Celine's technique is the
computation of recurrences for multiple sums $\sum_\kvec F(\nvec, \kvec)$, where
$F(\nvec, \kvec)$ is \qhg in all of its arguments.
For this we proceed by computing a so-called \kfree recurrence for the summand
first.

\begin{definition}
  A \qhg function $F(\nvec,\kvec)$ satisfies a
  \textit{\kfree recurrence}, if there exist a finite set $S$ of integer
  tuples of length $s+r$ and polynomials $\sigma_{\ivec,\jvec}(\nvec)$
  not all zero, such that
  \begin{equation} \label{recsumm}
     \sum_{(\ivec,\jvec) \in S} \sigma_{\ivec,\jvec}(\nvec) \, F(\nvec-\ivec, \kvec-\jvec) = 0
  \end{equation}
  holds at every point $(\nvec,\kvec)$ where all values of $F$ occurring in
  \eqref{recsumm} are well-defined.
  The set $S$ is called a \textit{structure set}.
\end{definition}

The computation of a \kfree recurrence is done by making an Ansatz of the form
\eqref{recsumm} for some structure set $S$ and undetermined $\sigma_{\ivec,\jvec}$.
Dividing equation~\eqref{recsumm} by $F(\nvec,\kvec)$ we get the rational equation
\begin{equation} \label{recsummrat}
  \sum_{(\ivec,\jvec) \in S} \sigma_{\ivec,\jvec}(\nvec) \,
  R_{F,\ivec,\jvec}(\nvec, \kvec) = 0,
\end{equation}
which after clearing denominators turns into the polynomial equation
\begin{equation} \label{recsummpoly}
   \sum_{(\ivec,\jvec) \in S} \sigma_{\ivec,\jvec}(\nvec) \, P_{F,\ivec,\jvec}(\nvec, \kvec) = 0.
\end{equation}
Next we compare the coefficients of all power products $q^{k_1 l_1}
\cdots q^{k_r l_r}$ in~\eqref{recsummpoly} with zero to get a homogeneous
system of linear equations for the $\sigma_{\ivec,\jvec}(\nvec)$.
Note that in this system the coefficients are rational functions in several
variables and not simply numbers.
Every non-trivial solution of the system gives rise to a \kfree recurrence.

\medskip
It has been shown by Wilf and Zeilberger~\cite{WilfZeilberger:AlgoProof} that
every so-called $q$-proper hypergeometric function $F(\nvec, \kvec)$ satisfies
a \kfree recurrence over some ``rectangular" structure set
$S_{\Ivec,\Jvec} := \{(\ivec,\jvec) \mid 0 \le i_l \le I_l,
0 \le j_m \le J_m\}$, since for large enough $\Ivec$ and $\Jvec$ the
number of unknowns exceeds the number of equations.

\medskip
However, this result is important only from theoretical point of view, because
in practice the run-time and memory demand grow very fast with the size of the
structure set $S$.
In particular, it turns out that rectangular structure sets are in general not usable,
since most recurrences live over a different domain, i.e.\ many
$(\ivec, \jvec) \in S_{\Ivec, \Jvec}$ are superfluous points with $\sigma_{\ivec,\jvec} = 0$.

Hence we also generalized the concept of \textit{$P$-maximal} structure sets
to the $q$-case, leading to more satisfactory results.
The underlying existence theory was originally introduced by
Verbaeten~\cite{Verbaeten:Diss}
for single-sums in the $q=1$ case.
Since it is based on arguments from plane geometry, there is no direct
generalization to multi-sums.
Nevertheless, as Wegschaider~\cite{Wegschaider:Dipl} pointed out,
$P$-maximal structure sets can be computed also in this situation:
the idea is to start with a small rectangular structure set and then
to add all those points that do not increase the degree of the polynomial
on the left-hand side of~\eqref{recsummpoly}.
This way the number of equations in the corresponding linear system remains
the same, whereas we maximize the number of unknowns.

However, there are still some cases where also this method, called
\textit{Verbaeten completion}, misses the minimal structure set.
In particular, this happens with most of the identities from Section~\ref{sec:triple}.
Hence we omit the details (see Riese~\cite{R}) and
only remark that in these specific instances we could overcome the problem
by first computing the structure set for $q=1$ and using the same for the
$q$-case.
Right now we do not have an explanation why this actually works.

\medskip
We want to emphasize that the correctness of a \kfree recurrence computed
by our program can be checked independently.
For that one simply divides equation~\eqref{recsumm} by $F(\nvec,\kvec)$
and verifies the resulting rational function identity~\eqref{recsummrat}
formally.
As Wegschaider showed this implies the correctness of the \kfree
recurrence even at those points $(\nvec, \kvec$) where $F(\nvec, \kvec) = 0$.

\medskip
Once we have computed a \kfree recurrence for the summand, the recurrence for
the whole sum can be obtained by summing over it.
Doing so the left-hand side of the recurrence could collapse to $0$.
Since we do not know of any example where this actually occurs, we do not go
further into the details here (see Wegschaider~\cite{Wegschaider:Dipl} or
Riese~\cite{R}).

Moreover, one should keep in mind that for some specific $\nvec$ certain
coefficients in the recurrence for the sum might vanish.
In this case one possibly has to consider extra boundaries.

\section{A Computer Proof of the Triple Bounded Identity} \label{sec:triple}

In this section we prove algorithmically identity~\eqn{1.15}
which we restate for convenience as
\begin{equation}\label{g=p}
  g_{i,j,k}(L_1,L_2,M) = p_{i,j,k}(L_1,L_2,M).
\end{equation}
For this we proceed in two steps.
First we compute a recurrence for $g_{i,j,k}(L_1,L_2,M)$ and show
that $p_{i,j,k}(L_1,L_2,M)$ satisfies the same recurrence.
Then we prove that \eqref{g=p} holds at a certain boundary.
The latter is achieved by repeating both steps for the boundary identity.

\subsection{The Recurrence}

First of all we load the package:
\small\begin{alltt}
\In{1} <<qMultiSum.m
\end{alltt}
\begin{alltt}
\Out{1} Axel Riese's qMultiSum implementation version 2.1 loaded
\end{alltt}\normalsize

\noindent
Then we enter the constraints and the summands.
\small\begin{alltt}
\In{2} a = i-ab-ac; b = j-ab-bc; c = k-ac-bc;
       t = a+b+c+ab+ac+bc;
       T[m_] := m(m+1)/2;
\end{alltt}
\begin{alltt}
\In{3} gsum = q^(T[t]+T[ab]+T[ac]+T[bc-1]) qBinomial[L\mysub{1}-t+a,a,q] *
              qBinomial[L\mysub{2}-t+b,b,q] qBinomial[M-t+c,c,q] *
              qBinomial[L\mysub{2}-t,ab,q] qBinomial[M-t,ac,q] qBinomial[M-t,bc,q] *
              (q^bc + (1-q^a)/(1-q^(L\mysub{1}-t+a)) (1-q^bc)/(1-q^(M-t-bc+1)));
\end{alltt}
\begin{alltt}
\In{4} psum = q^(s(M+2)-T[s]+T[i-s]+T[j-s]+T[k-s]) *
              qBinomial[L\mysub{1}-s,i-s,q] qBinomial[L\mysub{2}-i,j-s,q] *
              qBinomial[L\mysub{2}-i-j+s,s,q] qBinomial[M-i-j,k-s,q];
\end{alltt}
\normalsize

\noindent
The function for computing \kfree recurrences is called \texttt{qFindRecurrence}
(or \texttt{qFR} for short).
It takes as arguments the summand, the recurrence variables, the summation
variables, the dimensions of the initial rectangular structure set over
which Verbaeten completion is performed, and some optional parameters.
In our case we specify a structure set explicitly, which ---
as mentioned above --- we obtained by investigating the $q=1$ case.

\small\begin{alltt}
\In{5} qFindRecurrence[gsum, {\lbr}L\mysub{1},L\mysub{2},M,i,j\rbr, {\lbr}ab,ac,bc\rbr,
                       {\lbr}0,0,0,0,0\rbr, {\lbr}0,0,0\rbr,
                       StructSet -> {\lbr}{\lbr}0,0,0,0,0,0,0,0\rbr, {\lbr}0,1,0,0,0,0,0,0\rbr,
                                     {\lbr}1,1,1,0,1,0,0,0\rbr, {\lbr}1,2,1,1,1,1,0,0\rbr\rbr]
\end{alltt}
\Out{5}
\begin{align*}
  & q^{k + L_2}\,F(-1 + L_1,-2 + L_2,-1 + M,-1 + i,-1 + j,-1 + ab,ac,bc) + {} \\
  & q^{k + L_2}\,F(-1 + L_1,-1 + L_2,-1 + M,i,-1 + j,ab,ac,bc) + {} \\
  & q^k\,F(L_1,-1 + L_2,M,i,j,ab,ac,bc) - q^k\,F(L_1,L_2,M,i,j,ab,ac,bc) = 0
\end{align*}
\normalsize

\noindent
Then we sum up this recurrence by applying the function \texttt{qSumRecurrence}
(or \texttt{qSR} for short) to the previous result.
Of course we have to provide the information that the first $5$ variables are
recurrence variables.
Since we prefer backward shifts, we also call the function \texttt{BackwardShifts}.

\small
\begin{alltt}
\In{6} qSumRecurrence[%, 5] // BackwardShifts
\end{alltt}
\Out{6}
\begin{myhack}
\begin{equation}\begin{split}\label{grec}
   & q^{L_2}\,\SUM(-1+L_1,-2+L_2,-1+M,-1+i,-1+j) + {}\\
   & q^{L_2}\,\SUM(-1+L_1,-1+L_2,-1+M,i,-1+j) + {} \\
   & \SUM(L_1,-1 + L_2,M,i,j) - \SUM(L_1,L_2,M,i,j) = 0
\end{split}\end{equation}
\end{myhack}
\normalsize

Next we want to show that also $p_{i,j,k}(L_1,L_2,M)$ satisfies
the same recurrence.
Surprisingly, with our package we are not able to find~\eqref{grec}
directly.
On the other hand verification is a trivial task:
we simply show that already the summand of $p$ fulfills~\eqref{grec}
and therefore also $p$ itself.
For this we call the function \texttt{qCheckRecurrence} (or
\texttt{qCR} for short).
\small
\begin{alltt}
\In{7} qCheckRecurrence[%, psum]
\end{alltt}
\begin{alltt}
\Out{7} True
\end{alltt}
\normalsize

\subsection{The Boundary} \label{sec:Boundary}

First we note that identity~\eqref{g=p} is true
if any of the parameters $i,j,k$ is negative, since both
sides vanish in this case.
Thus, to complete the proof it suffices to show that
\begin{equation}\label{g=p bound}
  \tilde{g}_{i,j,k}(L_1,M) = \tilde{p}_{i,j,k}(L_1,M),
\end{equation}
where
\[
  \tilde{g}_{i,j,k}(L_1,M) := g_{i,j,k}(L_1,i+j-1,M)
\]
and
\[
  \tilde{p}_{i,j,k}(L_1,M) := p_{i,j,k}(L_1,i+j-1,M).
\]

To see this we have to distinguish two cases. If
$L_2 > i+j-1$ we rewrite recurrence~\eqref{grec}
as
\begin{align*}
  & \SUM(L_1,L_2,M,i,j) \\
  & \quad = q^{L_2}\,\SUM(-1+L_1,-2+L_2,-1+M,-1+i,-1+j) + \dotsb
\end{align*}
and find that on the new right-hand side either $j$ or $L_2$ is
shifted backwards by~$1$.
Analogously, if $L_2 \le i+j-1$ we rewrite \eqref{grec}
as
\begin{align*}
  & \SUM(L_1,-1+L_2,M,i,j) \\
  & \quad = {-q^{L_2}\,\SUM(-1+L_1,-2+L_2,-1+M,-1+i,-1+j)} - \dotsb
\end{align*}
and find that on the new right-hand side either $j$ is shifted
backwards or $L_2$ is shifted forwards by~$1$.
Hence, if we show \eqref{g=p bound} for all integer
valued $L_1$, then identity~\eqref{g=p} holds for all
integer valued $L_1$ and $L_2$.

\medskip
To prove \eqref{g=p bound} we proceed as before.
This means that we compute a recurrence for $\tilde{g}_{i,j,k}(L_1,M)$, check that
$\tilde{p}_{i,j,k}(L_1,M)$ satisfies the same recurrence and
prove \eqref{g=p bound} at a certain boundary.

\subsubsection{The Recurrence for the Boundary}

This time we succeed without spying out the structure set from
the $q=1$ case:

\small
\begin{alltt}
\In{8} qFR[gsum /. L\mysub{2} -> i+j-1, {\lbr}L\mysub{1},M,i\rbr, {\lbr}ab,ac,bc\rbr, {\lbr}0,0,0\rbr, {\lbr}0,0,1\rbr] //
         qSR[#, 3]& // BackwardShifts
\end{alltt}
\Out{8}
\begin{myhack}
\begin{equation} \label{gboundrec}
   q^{L_1}\,\SUM(-1 + L_1,-1 + M,-1 + i) +
   \SUM(-1 + L_1,M,i) - \SUM(L_1,M,i) = 0
\end{equation}
\end{myhack}

\begin{alltt}
\In{9} qCR[%, psum /. L\mysub{2} -> i+j-1]
\end{alltt}
\begin{alltt}
\Out{9} True
\end{alltt}
\normalsize

\subsubsection{The Boundary of the Boundary}

Finally, to complete the proof of \eqref{g=p bound} and thus of
\eqref{g=p}, we prove the identity
\begin{equation} \label{g=p boundbound}
  \tilde{g}_{i,j,k}(i+j-1,M) = \tilde{p}_{i,j,k}(i+j-1,M).
\end{equation}
Again we see that this is sufficient by viewing
recurrence~\eqref{gboundrec} as
\[
   \SUM(L_1,M,i) = q^{L_1}\,\SUM(-1 + L_1,-1 + M,-1 + i) +
   \SUM(-1 + L_1,M,i)
\]
in case of $L_1>i+j-1$, and as
\[
   \SUM(-1 + L_1,M,i) = -q^{L_1}\,\SUM(-1 + L_1,-1 + M,-1 + i) +
   \SUM(L_1,M,i)
\]
if $L_1 \le i+j-1$.
Hence, if we prove \eqref{g=p boundbound}, then identity~\eqref{g=p bound}
holds for all integer valued $L_1$.

\medskip
{}From the single-sum we immediately read off the relation
\[
  \tilde{p}_{i,j,k}(i+j-1,M) = \delta_{i,0} \, \delta_{j,0} \,
  q^{T_k} \qBin{\Delta}{k}{q},
\]
where $\Delta = M-i-j$.
Hence our boundary identity to verify becomes
\begin{equation}\label{g=p bound2}
  \tilde{g}_{i,j,k}(i+j-1,M) = \delta_{i,0} \, \delta_{j,0} \,
  q^{T_k} \qBin{\Delta}{k}{q}.
\end{equation}
The reason for switching from $M$ to $\Delta$ here is that
proving identity~\eqref{g=p bound2} for $\Delta=0$ is
easy, whereas for $M=0$ it is not at all.

\small
\begin{alltt}
\In{10} qFR[gsum /. {\lbr}L\mysub{1} -> i+j-1, L\mysub{2} -> i+j-1, M -> \(\Delta\)+i+j\rbr,
            {\lbr}\(\Delta\),i,j,k\rbr, {\lbr}ab,ac,bc\rbr, {\lbr}0,0,0,0\rbr, {\lbr}0,0,0\rbr,
            StructSet -> {\lbr}{\lbr}0,0,0,0,0,0,0\rbr, {\lbr}0,1,1,1,1,0,0\rbr, {\lbr}1,0,0,0,0,0,0\rbr,
                          {\lbr}1,0,0,1,0,0,0\rbr, {\lbr}1,0,1,1,0,0,0\rbr, {\lbr}1,0,1,1,0,0,1\rbr,
                          {\lbr}1,1,0,1,0,0,0\rbr, {\lbr}1,1,0,1,0,1,0\rbr, {\lbr}1,1,1,1,0,0,0\rbr,
                          {\lbr}1,1,1,2,0,1,1\rbr\rbr] //
           qSR[#, 4]& // BackwardShifts
\end{alltt}
\Out{10}
\begin{myhack}
\begin{equation}\begin{split} \label{grec bound2}
   & q^{-1 + 2\Delta + 2i + 2j}\,\SUM(-1 + \Delta,-1 + i,-1 + j,-2 + k) + {}\\
   & q^{-3 + \Delta + 3i + 3j}\,\SUM(-1 + \Delta,-1 + i,-1 + j,-1 + k) - {} \\
   & q^{\Delta + i + j}\,( -1 + q^{-1 + i + j}) \,\SUM(-1 + \Delta,-1 + i,j,-1 + k) - {} \\
   & q^{\Delta + i + j}\,( -1 + q^{-1 + i + j}) \,\SUM(-1 + \Delta,i,-1 + j,-1 + k) + {} \\
   & q^{\Delta + i + j}\,\SUM(-1 + \Delta,i,j,-1 + k) + \SUM(-1 + \Delta,i,j,k) - {} \\
   & q^{-1 + \Delta + 2\,i + 2\,j}\,\SUM(\Delta,-1 + i,-1 + j,-1 + k) -
   \SUM(\Delta,i,j,k) = 0
\end{split}
\end{equation}
\end{myhack}
\normalsize

\noindent
Obviously $\delta_{i,0} \, \delta_{j,0} \, q^{T_k} \qBinsm{\Delta}{k}{q}$ is
a solution of this recurrence, we only have to consider the two non-degenerate
cases when $i=j=0$ or $i=j=1$.
Note that once we have proved the validity of~\eqref{g=p bound2} for $\Delta=0$,
which follows immediately from~\cite[(3.7)]{AB1},
our recurrence implies the validity both for $\Delta > 0$ and
$\Delta \le 0$.
Indeed, for the case $\Delta > 0$ we rewrite \eqref{grec bound2} as
\[
  \SUM(\Delta,i,j,k) =
  q^{-1 + 2\Delta + 2i + 2j}\,\SUM(-1 + \Delta,-1 + i,-1 + j,-2 + k) + \dotsb
\]
and observe that on the new right-hand side at least one of
the parameters among $i, j, k$ is shifted backwards or
$\Delta$ is shifted backwards by $1$.
Analogously, for the case $\Delta \le 0$ we rewrite \eqref{grec bound2} as
\[
  \SUM(-1 + \Delta,i,j,k) =
   - q^{-1 + 2\Delta + 2i + 2j}\,\SUM(-1 + \Delta,-1 + i,-1 + j,-2 + k) - \dotsb
\]
and observe that on the new right-hand side at least one of the
parameters among $i, j, k$ is shifted backwards or $\Delta$ is shifted
forwards by $1$.

\medskip
Finally, let us briefly summarize the run-times.
The computation of~\eqref{grec bound2} took approximately $35$ seconds on an
\textsf{SGI Octane} using \textsf{Mathematica 4.0.1}.
All other computations could be carried out within less than $5$
seconds.

\section{Computer Proofs of New Finite Analogs\\ of Jacobi and Euler Formulas}
\label{sec:Euler}

In~\cite[(5.6)]{AB1} Alladi and Berkovich pointed out that \eqn{1.14} can be 
used to derive
\begin{equation} \label{finiteJac}
   \sum_{l=0}^L a^{-l} \; \frac{1+a^{2l+1}}{1+a} \; q^{T_l} =
   \sum_{i,j,k} a^{i-j} \, (-1)^k \, q^{T_i+T_j+T_k} \qBin{L-k}{i}{q}
   \qBin{L-i}{j}{q} \qBin{L-j}{k}{q},
\end{equation}
from which Jacobi's formula
\[
  \sum_{l \ge 0} (-1)^l \, (2l+1) \, q^{T_l} = (q;q)_\infty^3,
\]
follows in the limit $L \to \infty, a \to -1$.
It is instructive to prove \eqref{finiteJac} in an automated fashion.
For the right-hand side of \eqref{finiteJac} we obtain within $20$ seconds the
following recurrence of order $4$:
\small
\begin{alltt}
\In{11} Clear[a]
\end{alltt}
\begin{alltt}
\In{12} qFR[a^(i-j) (-1)^k q^(T[i]+T[j]+T[k]) qBinomial[L-k,i,q] *
            qBinomial[L-i,j,q] qBinomial[L-j,k,q],
            L, {\lbr}i,j,k\rbr, 2, {\lbr}0,0,0\rbr] // qSR
\end{alltt}
\Out{12}
\begin{align*}
  & {-\big( a\,q^{9 + 3L}\,\SUM(L) \big)}  +
   q^{7 + 2L}\,( 1 - a + a^2 + a\,q^{2 + L}) \,\SUM(1 + L) -  {} \\
  & ( 1 - a + a^2) \,q^{4 + L}\,( -1 + q^{3 + L}) \,\SUM(2 + L) + {} \\
  & ( -a - q^{4 + L} + a\,q^{4 + L} - a^2\,q^{4 + L})\,\SUM(3 + L) +
    a\,\SUM(4 + L) = 0
\end{align*}
\normalsize
Now we plug in the left-hand side of \eqref{finiteJac}:
\small
\begin{alltt}
\In{13} Simplify[% /. SUM[L + m_] :> SUM[L] +
                 Sum[a^(-l) (1+a^(2l+1))/(1+a) q^T[l], {\lbr}l,L+1,L+m\rbr]]
\end{alltt}
\begin{alltt}
\Out{13} True
\end{alltt}
\normalsize
Once we have verified \eqref{finiteJac} for $L \in \{0,1,2,3\}$, we are done.

\medskip
Recently, we came up with the identity
\begin{equation} \label{finiteEuler}
   \sum_{l=0}^L q^{2(T_L-T_l)} =
   \sum_{i,j,k} (-1)^j \, q^{2T_i+2T_j+2T_k-i-j} \qBin{L-k}{i}{q^2}
   \qBin{L-i}{j}{q^2} \qBin{L-j}{k}{q^2},
\end{equation}
which for $L \to \infty$ turns into Euler's formula
\[
  1 = (-q;q)_\infty \, (q;q^2)_\infty.
\]
With our package the proof can be done again in $20$ seconds:
\small
\begin{alltt}
\In{14} qFR[(-1)^j q^(2T[i]+2T[j]+2T[k]-i-j) qBinomial[L-k,i,q^2] *
            qBinomial[L-i,j,q^2] qBinomial[L-j,k,q^2],
            L, {\lbr}i,j,k\rbr, 2, {\lbr}0,0,0\rbr] // qSR
\end{alltt}
\Out{14}
\begin{align*}
   & {-\big( q^{14 + 6L}\, \SUM(L) \big)}  +
   q^{12 + 4L}\, ( 1 + q^{4 + 2L} ) \,\SUM(1 + L) - {} \\
   & q^{6 + 2L}\, ( -1 + q^{3 + L} ) \, ( 1 + q^{3 + L}) \,
   \SUM(2 + L) + {}\\
   & ( -1 - q^{8 + 2L} ) \,\SUM(3 + L) +
   \SUM(4 + L) = 0
\end{align*}
\normalsize
Plugging in the left-hand side of~\eqref{finiteEuler} gives:
\small
\begin{alltt}
\In{15} Simplify[% /. SUM[L + m_] :> q^(2(L m + T[m])) *
                 (SUM[L] + Sum[q^(2(T[L]-T[l])), {\lbr}l,L+1,L+m\rbr])]
\end{alltt}
\begin{alltt}
\Out{15} True
\end{alltt}
\normalsize
Again the remaining task is to verify \eqref{finiteEuler} for $L \in \{0,1,2,3\}$.

\section{Concluding Remarks} \label{sec:Conclusion}

While it is apparent that the package \texttt{qMultiSum} is a valuable tool
for proving $q$-hypergeometric identities, we would like to  
point out that our proof of~\eqn{1.15} is not fully automated.
First of all, in general, the computer generates not just one, but many recurrences.
Some of these recurrences are  
``dead ends" in practical terms.

For instance, if one starts with the simple recursion relation
\begin{equation} \label{eq:5.1}
   g_{i,j,k}(L_1,L_2,M)=g_{i,j,k}(L_1-1,L_2,M)+q^{L_1}g_{i-1,j,k}(L_1-1,L_2-1,M-1) 
\end{equation}
and then sets $L_1=i-1$, one gets
\begin{equation} \label{eq:5.2}
   g_{i,j,k}(i-1,L_2,M)=q^{i(M+2)-T_i+T_{j-1}+T_{k-i}}
   \qBin{L_2-i}{j-i}{q} \qBin{L_2-j}{i}{q} \qBin{M-i-j}{k-i}{q}.
\end{equation}
Unfortunately, we did not succeed in proving the boundary identity~\eqn{5.2} with 
\texttt{qMultiSum} and, as a result, we were not able to complete the proof
of~\eqn{1.15}, taking~\eqn{5.1} as the starting point.

Second, crucial ``human insight" was used in selecting the not so obvious boundaries
\[ 
  L_1=i+j-1, \quad L_2=i+j-1, \quad M=i+j
\]
in Section~\ref{sec:Boundary}.
These boundaries were determined by two requirements: 
\begin{enumerate}
\item The validity of the identity on the boundary, along with the recursion relations,
      should imply the validity of the identity everywhere.
\item The identity should take a particularly simple form on the chosen boundary.
\end{enumerate}
We have already mentioned at the end of Section~\ref{sec:qMultiSum} that the recursion
relations generated by the computer should be examined carefully, because on some
specific hyperplanes these recurrences may become a triviality $0=0$.
At present, this examination has to be carried out manually.

\medskip
We would like to finish this article by pointing out that currently only \eqn{1.14}
finds a partition theoretical interpretation.
Clearly, it is highly desirable to find a partition theorem which corresponds to
identity~\eqn{1.15}.
It is encouraging that for the case when one of the parameters $i,j,k$ in~\eqn{1.15}
is set to zero, such a theorem was recently found in~\cite{AB2}.

\subsection*{Acknowledgment}

We would like to thank Krishnaswami Alladi, George E.~Andrews and Doron Zeilberger for 
their interest and comments on the manuscript.

\end{document}